\documentclass[11pt,a4paper,final,leqno,notitlepage]{amsart}
\usepackage[dvips,final]{epsfig}
\usepackage{epsf, graphics, latexsym, eucal, pstcol,
amsfonts, amstext, amssymb, amsmath, amsthm, amsbsy, amsxtra, amscd }

\allowdisplaybreaks[1]

\theoremstyle{plain}

\newtheorem{ltheorem}{Theorem}           
\newtheorem*{corollary}{Corollary}    

\newtheorem{lproposition}{Proposition}

\theoremstyle{definition}

\newtheorem{ldefinition}{Definition}

\theoremstyle{remark}

\newtheorem*{remark}{Remark}

\begin{document}

\hoffset = 0pt
\voffset = 0pt
\topmargin = 0pt
\headheight = 0pt
\textwidth = 355pt

\setlength{\leftmargin}{0.1cm}
\setlength{\rightmargin}{0.1cm}

\renewcommand{\abstractname}{Abstract}

\title{Dual scattering channel schemes and transmission line methods}

\author{Steffen Hein}

\keywords{Time domain methods, finite difference schemes, scattering,
hyperbolic PDE, TLM, FDTD, DSC. \hfill MSC-classes: 65L12, 65L20, 65M06}

\begin{abstract}

Dual scattering channel (DSC) schemes naturally extend the transmission
line matrix (TLM) numerical method beyond the lines set up by P.B. Johns
and coworkers.
Conceptually, DSC schemes retain from TLM the typical formal splitting of
the computed fields into incident and reflected components, interchanged
and scattered between adjacent cells, without yet a need for transmission
lines.
They generally admit nontrivial cell interface scattering and so release
from a set of modeling limitations to that the classical TLM method is
subject.
The elements characterizing DSC schemes are outlined and illustrated by
coupling within a non-orthogonal mesh a model of heat propagation to a lossy
Maxwell field. 

\vspace{.2cm}
\hspace{-.50cm}
Westerham on \today
\end{abstract}
\maketitle
\pagestyle{myheadings}
\markboth{{\normalsize \textsc{Steffen Hein}}}
{{\normalsize \textsc{Dual Scattering Channel Schemes}}}
\normalsize
\vspace{-0.6cm}
\begin{quote}\notag 
\small
{
\textnormal{Life is not an illogicality;
yet it is a trap for logicians.} \newline
\textnormal{It looks just a little more mathematical
and regular than it is.}
\begin{flushright}
\textit{Gilbert Keith Chesterton}
\end{flushright}
}
\end{quote}
\vspace{-0.2cm}
\section{Introduction}\label{S:int}
Chesterton's trap lurs where theory encouters real life - ~or other theories.
Mathe\-mati\-cians apply their methods, and sometimes export them to fields
where they are inadequate; just as physicists or microwave engineers
introduce proper concepts into mathematics, with sometimes creative effects
or walking into a trap - ~two issues that are not absent in the history
of the transmission line matrix (TLM) numerical method.

In the early 1970s Peter B. Johns\,\dag, prominent professor of electrical
engineering at the University of Nottingham, had with his coworkers the
brilliant idea to apply transmission line network modeling techniques to
the numerical solution of Maxwell's equations ~\cite{JoBe, Jo1}.
There are in fact appealing features of the TLM method,
such as unconditional stability ~\cite{Jo2}, high flexibility due to
novel stub-loading (deflection) techniques ~\cite{Ch,Hoe,He1}, and easy
implementation, which made the method attractive to users.
TLM approximation to wave propagation, diffusion, and transport phenomena
became so a subject of assiduous study and publication [Tlm1-3].

Familiarity with the transmission line picture and the well known
scattering concept certainly fostered the acceptance of the TLM method
in the microwave engineering community.
However, the primary interest so turned of course on applications in their
own discipline, rather than onto the inner structure of the TLM algorithm
as an object of mathematical analysis.
The underlying framework remained firmly linked to 'real' transmission
lines, indeed too firmly to allow for substantial evolution of the TLM
method beyond the original field of application.
The algorithm could not unfold its entire potentiality.

The history of the TLM method is best displayed in the proceedings of
three conferences that exclusively focussed on the subject [Tlm1-3].
Also, the monographs of Christopoulos ~\cite{Ch} and de ~Cogan ~\cite{dC}
reproduce the original ideas and introduce into classical applications.
Rebel ~\cite{Re} gives a fairly complete survey over the state of the
art by the year 2000 - ~which is more or less still the actual state.

Dual scattering channel (DSC) schemes result from an incisive revision
of the TLM paradigm.
They are generalized transmission line methods in arising from a twofold
abstraction:
Firstly, the scattering channel concept underlying TLM is redefined in
terms of paired distributions. 
(Characteristic impedances are thus neither needed, nor in general
defined, e.g.)
In the second place, nontrivial cell interface scattering is admitted
during the connection step of iteration thus taking advantage of the
intrinsic duality in the connection-reflection cycle of the algorithm.

The new framework bypasses a set of modeling limitations imposed by
transmission lines while it preserves the major advantages of the TLM
method. In particular, the convolution type updating scheme and
Johns' two step connection-reflection cycle are essentially retained.
DSC schemes remain unconditionally stable under quite general
circumstances, and they allow for process parallelization and deflection
- ~the latter of which generalizes the stub loading feature of the classical
TLM method ~\cite{He1, He2}.
 
DSC schemes are conceptually simple and natural. They are a coherent
response to space time discretization in that meshing virtually
\emph{imposes\/} an artificial cell-boundary duality.
A detailed technical introduction into DSC schemes has recently been
given in ~\cite{He3}. The present paper offers a condensed desription of
the framework with emphasis on the typical contraction properties that
ensure stability.

\textsc{Spinner}'s implementation of a DSC heat propagation scheme in 
non\--or\-tho\-go\-nal mesh, wherein heat sources are coupled to a
lossy Maxwell field, illustrates the approach.
\vspace{-0.2cm}
\section{Causality}\label{S:sec0}

Let ${\mathcal{L}\/}$ be a real or complex linear space and $I$ 
a totally ordered set 
(\,e.g. $I \in \{\mathbb{N}\,, \mathbb{Z}\,, \mathbb{R}\}\,$;
intervals are then naturally defined in $\,I\,$ by means of the order relation.
Of course, we shall think with $\,I\,$ of a discrete or continuous time domain.)
Also, let ${\mathcal{E} \subset \mathcal{L}^{I}\/}$ be a set
of functions such that ${f\in\mathcal{E}\,}$ implies
$\chi_{_{s\leqslant t}}(s) f(s)\in\mathcal{E}\,$,
for every ${t\in I\,}$\,, where $\chi_{_P}(s)$ denotes the
characteristic function of property $\,P\,$\, 
(which is 1 if $\,s\,$ shares that property and 0 else).
Sometimes we write simply \quad $[\,f\,]_{_{\leq t}}$ \quad
in the place of \quad $\chi_{_{s\leqslant t}}(s)\;f(s)$.

\begin{ldefinition}\label{D.0.1}
A function
$F:\mathcal{E}\to\mathcal{E}\,$ is called \emph{causal\/},
iff for every ${f\in\mathcal{E}\,}$\\ and ${t \in I\,}$
\vspace{-.25cm}
\begin{equation}\centering\notag 
F\,f\,(\,t\,)\;=\;F\;[\;\chi_{_{s\leqslant t}}f\,(\,s\,)\;]\,(\,t\,)\quad .
\end{equation}
Such functions are also called (causal) \emph{propagators}.

\newpage
\begin{remark}
In some respect, causal functions generalize lower triangular matrices
or integral operators such as
\vspace{-.15cm}
\begin{equation}\centering\notag 
F\,f\,(\,t\,) \; = \; \int_{-\infty}^{t}K\,(\,t\,-\,s\,)\,f\,(\,s\,)\; ds\;,
\quad
\end{equation}
for instance with a Green's function kernel $K$. \\
Note that in general $\,F:f\mapsto F f\,$ needs not to be linear.

Typically $\,f\,$ represents a state evolving in time
(\,i.e. a \emph{process}\,).
Then causality of $\,F\,$ spells that $\,F\,f\,(\,t\,)\,$
depends on the history of $\,f\,$ only up to present 
time $\,t\,$.
\end{remark}
\end{ldefinition}

\hspace{-.53cm}
The proof of the following is easy and left to the reader.
\begin{lproposition}\label{P.0.1}\ \\
For every ${\,t\in I\,}$;
${\,f,\,g\,\in\mathcal{E}\,}$ and causal functions
${F,\,G\,:\,\mathcal{E}\to\mathcal{E}}$
\begin{itemize}
\item[(i)]
$F f\,(\,s\,)\;=\;F\;[\;\chi_{_{r\leqslant t}}(\,r\,)\,f\,(\,r\,)
\;]\,(\,s\,)\,$,
for every $\,s\leqslant\,t$ \;.
\vspace{4pt}
\item[(ii)]
If $\;f\,(s)\;=\;g\,(s)\;$ for \emph{every} $\,s\leqslant\,t\;$, \\
\vspace{-0.45cm}
\begin{equation}\centering\notag 
\text{then} \quad F f\,(s)\;=\;F g\,(s)\;\,\text{for every\,}
\,s\leqslant\,t\;.
\end{equation}
\end{itemize}
\begin{itemize}
\item[(iii)]
The product of causal functions
\begin{equation}\centering\notag 
\begin{aligned}
FG \, : \, \mathcal{E}\;&\longrightarrow\; \mathcal{E}\\[-2pt]
f\;&\longmapsto\;F\,G\,f\,\;:\,=\;F\,[\,G\,[\,f\,]\,]
\end{aligned}
\end{equation}
is again causal. In fact, if $\mathcal{E}$ is a linear space,
then the causal functions over $\mathcal{E}$ form an algebra.
\end{itemize}
\end{lproposition}
\vspace{-0.2cm}
\section{DSC Processes}\label{S:sec1}

Just as the TLM algorithm, DSC schemes operate on a linear space
$\,\mathcal{P}\,$
of so-called \emph{propagating fields\,}, which is a product of
(real or complex) normed spaces
\begin{equation}\centering\notag 
\mathcal{P}\;=\;\;\mathcal{P}_{in}\,\times\,\mathcal{P}_{out}\quad,
\end{equation} 
cf. ~\cite{He3}.
The two factors are (loosely) named the \emph{incident} and
\emph{outgoing} subspace of ${\,\mathcal{P}}$. 
They are isomorphic in that there exists a canonical involutary
isomorphism of normed spaces
\begin{equation}\centering\notag 
\begin{aligned}
nb\;:\;\mathcal{P}\;&\to\;\mathcal{P}\\[-2pt]
z\;=(\,z_{in}\,,\,z_{out}\,)\;&\mapsto\;
(\, z_{out}\,,\,z_{in}\,)\;=\,:\; nb\,(\,z\,)\;,
\end{aligned}
\end{equation}
which is commonly called the \emph{node-boundary} map. 
Hence, there exists a space
${(\,\mathcal{L}\,,\,\|...\|\,)\,}$
such that
\begin{equation}\centering\notag 
\mathcal{P}_{in}\;\cong\;\mathcal{P}_{out} 
\;\cong\;(\, \mathcal{L}\,, \,\|...\|\,)\quad \qquad
\end{equation}
in the sense of isomorpy of normed linear spaces. In the same meaning
\begin{equation}\centering\notag 
\begin{aligned}
\mathcal{P}\;\cong\;(\,\mathcal{L}^{2}\,, 
\, \|...\|^{\sim}\,) &\quad,\\[-2pt] 
\text{e.g. with norm} \quad
&\|\,(\,a\,,\,b\,)\,\|^{\sim}\;:\,=\;
\sqrt{\| a \|^{2}\,+\,\| b \|^{2}}
\quad ;\quad a,b \in\mathcal{L}\;
\end{aligned}
\end{equation}
(or any equivalent norm).

DSC and TLM algorithms then follow a two-step iteration cycle in working
with alternate application of a \emph{connection} and \emph{reflection} map
\vspace{-0.3cm}\newline
\makebox[6.5cm][s]
{
   \begin{minipage}[t]{6.5cm}
   \begin{equation}\centering\label{1.5}
   \begin{aligned}
   \mathcal{C}:I\times\mathcal{P}_{out}^{I}&\,\to\,\mathcal{P}_{in}\\
   (t,z_{out}(t-\kappa\tau)_{\kappa=0}^\infty)
   &\,\mapsto\,z_{in}(t)\;,
   \end{aligned}
   \end{equation}
   \end{minipage}
}
\;
\makebox[5.3cm][s]
{
   \begin{minipage}[t]{5.3cm}
   \begin{equation}\centering\notag
   \begin{aligned}
   \mathcal{R}:J\times\mathcal{P}_{in}^{J}&\,\to\,\mathcal{P}_{out}\\
   (t,z_{in}(t-\kappa\tau)_{\kappa=0}^\infty)
   &\,\mapsto\,z_{out}(t)
   \end{aligned}
   \end{equation}
   \end{minipage}
}
\vspace{0.2cm}\newline
which respectively update the propagating fields at even and odd
integer multiples of half the time step, i.e. on ${I\,:\,=\,}$ ${
\{\,k\tau\mid k\in\mathbb{Z}\,\}}$ and
${J\,:\,=\,}$ ${\{\,(2k+1)\tau/2\mid k\in\mathbb{Z}\,\}}$.
To these maps the following functions $\,F_{_{C}}\,,
F_{_{R}}\,$ are associated in a one-to-one correspondence
\begin{equation}\centering\notag 
\begin{aligned}
F_{_{C}}\,:\,\mathcal{L}^{I}\,&\to\,\mathcal{L}^{I} \\[-6pt]
f\,&\mapsto\,F_{_{C}}\,f\quad\text{with}\quad
F_{_{C}}\,f\,(\,t\,)\;:\,=\;nb\circ\mathcal{C}\,(\,t\,,\,f\,)
\quad
\end{aligned}
\end{equation}
and 
\begin{equation}\centering\notag 
\begin{aligned}
F_{_{R}}\,:\,\mathcal{L}^{J}\,&\to\,\mathcal{L}^{J} \\[-6pt]
g\,&\mapsto\,F_{_{R}}\,g\quad\text{with}\quad
F_{_{R}}\,g\,(\,t\,)\;:\,=\;nb\circ\mathcal{R}\,(\,t\,,\,g\,)
\quad.
\end{aligned}
\end{equation}

\begin{lproposition}\label{P.1.1}
\begin{itemize}
\item[]
\item[(i)]
$F_{_{C}}\,$ and $\,F_{_{R}}\,$ are causal on
$\mathcal{L}^{I}$ and $\mathcal{L}^{J}$, respectively.
\item[(ii)]
For every $\,r,s\in J\,$ and $\,T_{q}:f(t)\,\mapsto\,f(t+q)\,$
the \emph{shift operator} on $\,\mathcal{L}^{I\;\cup\,J}\,$,
\vspace{-6pt}
\begin{equation}\centering\notag 
\quad\text{if}\quad r+s\;\leqslant\;0\;,\quad\text{then}\;
\begin{cases}
T_{r}\circ F_{_{C}}\circ T_{s}\\[-3pt]
T_{r}\circ F_{_{R}}\circ T_{s}
\end{cases}
\text{is causal on}\quad
\begin{cases}
\mathcal{L}^{J}\\[-3pt]
\mathcal{L}^{I}
\end{cases}
\!.
\end{equation}
\end{itemize}
\end{lproposition}

\hspace{-.53cm}
The first statement is validated by observing that the reflection and
connection maps $\,\mathcal{R}$ and $\,\mathcal{C}$ are functions of back
in time running sequences of propagating fields, cf. ~\eqref{1.5}, and
are themselves causal in this sense ~\cite{He3}. The second statement is
then a directly verified. \hfill $\Box$

\vspace{0.2cm} 
\hspace{-.55cm} 
For later use we fix straightaway:
\begin{ldefinition}\label{D.1.1}
For any incident function
$\,e\,\in\,\mathcal{L}^{I}\times\{\,0\,\}\,$
supported on a finite interval $\,[\,0\,,\,N\tau\,)\subset I\,$,
the DSC process excited with $\,e\,$ and generated by
$\,\mathcal{R}\,$ and $\,\mathcal{C}\,$ is the unique function
$\,z\,=\,(\,z_{in}\,,\,z_{out}\,)\,\in\,(\,\mathcal{L}^{2}\,)^{I \cup J}\,$ 
such that $\,z\,(\,t\,)\,=\,0\,$ for $\,t\,\leqslant\,0\,$ and recursively
for $\,0\,<\,t\,\in\,I\cup J\,$
\vspace{-2pt}
\begin{equation}\centering\label{1.11}
\begin{aligned}
z\,(\,t\,+\,\frac{\tau}{2}\,)\,=\,
\begin{cases}
(\,z_{in}\,(\,t\,)\,,\;T_{-\frac{\tau}{2}}\,F_{_{R}}\,T_{-\frac{\tau}{2}}
\,[\,e\,+\,z_{in}\,]\,(\,t\,)\,)
&\,\text{if $\,t\,\in\,I\,$}\\
(\,T_{-\frac{\tau}{2}}\,F_{_{C}}\,T_{-\frac{\tau}{2}}\,
[\,z_{out}\,]\,(\,t\,)\,,\;z_{out}(\,t\,)\,)
&\,\text{if $\,t\,\in\,J\,$}\,.
\end{cases}
\end{aligned}
\end{equation}
\end{ldefinition}

Thus, in particular, $\,z_{in}$ and $\,z_{out}$ 'switch' respectively
at even and odd integer multiples of half the time step $\,\tau$.
\vspace{-0.2cm}
\section{Passivity}\label{S:sec2}

Algorithm stability prevents the computational process from piling up
to infinity (of course, it does not yet imply convergence,
or ~\emph{consistence} of the algorithm).
TLM models are \emph{unconditionally\/} stable in that they are
equivalent to passive linear transmission line networks ~\cite{Jo2}.
DSC schemes, in not using lines, need a more general characterization
which here is given in terms of $\alpha$-passive causal functions.

Under the same assumtions for $\,\mathcal{L}\,$ and $\,I\,$ as in 
section ~\ref{S:sec0}, let ${\|...\|\,}$ be a norm on $\mathcal{L}\,$ and 
$\alpha \in \mathbb{R}^{\,\mathcal{L}}\,$ a continuous non-negative real
functional.
\begin{ldefinition}\label{D.2.2}
\begin{itemize}
\item[]
\item[(i)]
We call a process ${g\,:\,I\to\mathcal{L}}$ \emph{stable},
iff it is uniformly bounded on $I\,$
(i.e.\! there exists ${\,b\in\mathbb{R}_{+}\,}$ such that 
${\|\,g\,(t)\,\|\,<\,b\,}$ for every ${t\in I\,}$\,).
\item[(ii)]
The functional ${\alpha\,:\,\mathcal{L}\to\,\mathbb{R}\,}$
is named a (\emph{de})\emph{limiting functional\,},
iff there exist any non-negative real constants 
$\,a,b,c\,$ such that
\begin{equation}\centering\label{2.4}
\|\,z\,\|\;\leqslant\;a\,+\,b\,(\alpha\,(\,z\,)\,)^{c}\quad,
\end{equation}
for every ${\,z\in\mathcal{L}\,}$.
Then obviously ${\,b, c > 0\,}$, and we say also that $\alpha\,$ 
is \emph{minimal increasing\,} (\emph{in any order\,})
\emph{not lesser than} ${\,1/c\,}$.
\end{itemize}
\end{ldefinition}

Let $\,\mu\,$ be a measure on $\,I\,$ such that intervals are
$\,\mu$-measurable sets.
Functions on $\,I\,$ are henceforth read \emph{modulo $\mu$}
(viz. as equivalence classes of functions that differ at most on sets of
$\mu$-measure zero). Also, let
${\,\alpha\in\mathbb{R}^{\,\mathcal{L}}\,}$
be a delimiting functional on $\mathcal{L}\,$
(\,increasing not lesser than any positive order ${\,1/c\,}$),
and assume that ${\,\alpha\circ f\,}$
is $\,\mu-$summable over finite intervals in $\,I\,$
for every ${\,f\in\mathcal{E}\,}$.
The latter is, for instance, the case if 
$\,\alpha\,(\,z\,)=\|\,z\,\|^{\,p}$
for any real ${\,p\geqslant 1\,}$
and ${\,\mathcal{E}\,\subset\,L^{\,p}(\,I,\,\mathcal{L}\,)\,}$
which is the metric completion of
\vspace{-3pt}
\begin{equation}\centering\notag 
\{\,f\,\in\mathcal{L}^{I}\,\mid\,
(\,\int\nolimits_{I}\|\,f\,\|^{p}\;d\mu\;)^{\;\frac{1}{p}}
\; \leqslant \; \infty \} \quad,
\end{equation}
i.e. $\mathcal{E}$ is a subset of the Banach space with norm
$\|\,f\,\|_{p}\;:\,=\;
(\,\int\nolimits_{I}\|\,f\,\|^{p}\;d\mu\,)^{\,\frac{1}{p}}$.

\begin{ldefinition}\label{D.2.3}
A causal function $\,F:\mathcal{E}\to\mathcal{E}\,$
is called $\alpha$-\emph{passive\/}, iff
\begin{equation}\centering\label{2.6}
\int\nolimits_{s<t}\alpha\,(\,F f\,(\,s\,)\,)\;d\mu(s)\;\leqslant\;
\int\nolimits_{s<t}\alpha\,(\,f\,(\,s\,)\,)\;d\mu(s)\quad,
\end{equation}
for every ${\,f\in\mathcal{E}\,}$ and ${\,t\in I\/}$.
\begin{remark}
If $\,\alpha\,=\,\|...\|^{p}$ for any real $\,p\geqslant 1\,$ and 
${\,\|\,f\,\|^{p}\,}$ is $\,\mu$-summable over $\,I\,$,
i.e. $\,f\in L^{\,p}(\,I,\,\mathcal{L}\,)\,$,
then \eqref{2.6} clearly implies 
$\,F f\, \in L^{\,p}(\,I,\,\mathcal{L}\,)\,$ and
\begin{equation}\centering\notag 
\|\,F f\,\|_{_{p}} \;\leqslant\;
\|\,f\,\|_{_{p}} \quad.
\end{equation}
Hence, every ${\|...\|^p}$-passive causal function $\,F\,$ 
defines a \emph{contraction} operator on
${\mathcal{E}\,\cap\,L^{\,p}(\,I,\,\mathcal{L}\,)}$.
\end{remark}
\end{ldefinition}

Assume furthermore that ${\tau\in\mathbb{R}_{+}\,}$ and let on
${\,I\;:\,=\;\mathbb{R}\,}$ the measure $\mu$ be concentrated in
${\{\,k\tau\mid k\in\mathbb{Z}\,\}\,}$ with uniform weight
${\mu\,(\,\{k\tau\}\,)\,=\,\tau\,}$, ${k\in\mathbb{Z}\,}$. 
Alternatively, let $I\;:\,=\;\{\,k\tau\mid k\in\mathbb{Z}\,\}\,$
with 'the same' measure $\,\mu\,$. 
(Virtually we deal of course with that discrete situation, 
even in working on the real axis with functions that are constant
over intervals
${[\,k\tau, (k+1)\,\tau)\,}$; ${k\in\mathbb{Z}\,}$.
We retain the integral formalism for simplicity, but the reader may
optionally re-write the following integrals as sums.)
For every $f\in\mathcal{E}$ let ${f(t-\tau)}$ $\in\mathcal{E}$,
i.e. $\,\mathcal{E}\,$ is closed under time shifts by negative
integer multiples of $\,\tau$.

For arbitrary ${N\in\mathbb{N}_{+}\,}$ and \emph{exciting function\,}
$e\in\mathcal{E}\subset\mathcal{L}^{I}\,$ with support on
${[\,0,\,N\tau\, )\subset I\,}$, the following holds

\begin{ltheorem}\label{T.2.1}
\emph{[\,Stability of the iterated passive causal process\,]}
\newline
For every $\alpha$-passive causal function
$\,F:\mathcal{E}\to\mathcal{E}\,$ and
${\,e \in \mathcal{E}\,}$ as stated,
\newline if ${g\in\mathcal{L}^{I}\,}$ is a process such that
${g(t)\,\equiv\,0\,}$ for $\,t\leqslant 0\,$ and
\begin{equation}\centering\label{2.8}
g\,(\,t\,+\tau\,)\;=\;F\,[\,e\,+\,g\,]\,(\,t\,)\quad
\end{equation}
for ${\,t\,=\,n\tau;\,n\in\mathbb{N}\,}$,
then $\,g\,$ is \emph{uniquely defined}
\textnormal{(}modulo $\mu$\textnormal{)},
and for every $t\geqslant N\tau\,$ holds
\vspace{-.25cm}
\begin{equation}\centering\label{2.9}
\|\,g\,(\,t\,)\,\|\;\;\leqslant\;\;a\;+\;
(\;\;\frac{b}{\tau}\;\int_{[\,0,\,N\tau\,)}\alpha\,(\,e\,+\,g\,) -
\alpha\,(\,g\,)\;\;\;d\mu\;\;)^{\;c}\quad\qquad
\end{equation}
with all constants $a,b,c\,$ that satisfy \eqref{2.4}.
Hence $g\,$ is \emph{stable}.
\end{ltheorem}

\begin{remark}
A process $\,g\in\mathcal{L}^{I}\,$ that is recursively generated
according to \eqref{2.8} by iteration of an $\,\alpha$-passive causal
function $\,F$ is called an $\alpha$-\emph{passive process}.

The theorem ensures thus that for any excitation of finite duration 
(and with no further restrictions) the $\alpha$-passive process
is necessarily stable and in this sense \emph{unconditionally} stable.

Note that \emph{existence} of such a process $g$ is not a priori guaranteed,
since \linebreak this obviously depends on the condition that with
${s_{0}\,:\,=\;e\,\in\mathcal{E}\,}$\linebreak
also the following functions belong to $\mathcal{E}\,$
\begin{equation}\centering\label{2.10}
s_{n}\;:\,=\;
[\,e\,+\,F[\,s_{n-1}\,]\,]_{_{\leq n\tau}}
\;\in\;\mathcal{E}\;;
\quad\text{for}\;0\,<\,n\,<\,N\;,
\end{equation}
which has to be checked if need be.
\newline Clearly, condition \eqref{2.10} is always true (hence $\,g\,$
exists) if $\,\mathcal{E}\,$ is a linear space.
We do not universally premise this, in order to apply the theorem also
to non-linear situations, where conditions \eqref{2.10} may only be
satisfied for sufficiently small excitations $e\,$.
\end{remark}

\begin{corollary}
\begin{itemize}
\item[]
\item[(i)]
In the special case $\,\alpha\,=\,\|...\|\,$ estimates \eqref{2.4} holds
with ${a=0}$, ${b=1}$, ${c=1}$.
Then the triangle inequality applies to the
integrand of \eqref{2.9} and validates the bound
\begin{equation}\centering\notag 
\,\|\,g\,(\,t\,)\,\| \;\; \leqslant \;\;
\frac{1}{\tau}\;\int_{[\,0,\,N\tau\,)}\|\,e\,\|
\;\;d\mu \quad.
\end{equation}
\item[(ii)]
If $\,N=1\,$, i.e. ${\,e(t)\,}$ is a \textsc{Dirac} excitation concentrated
on ${\,[\,0\,,\tau\,)\,}$ \,\textnormal{(}where ${\,g=0\,}$\textnormal{)},
then \eqref{2.9} reads simply
\begin{equation}\centering\notag 
\|\,g\,(\,t\,)\,\|\;\leqslant\;a\,+\,
(\;\frac{b}{\tau}\;\int_{[\,0,\tau\,)}\alpha\,(\,e\,)
\;\;d\mu\;\;)^{\,c} \quad,
\end{equation}
provided that $\,\alpha\,(\,0\,)\,=\,0\,$ \textnormal{(}
which is the normal case \textnormal{)}.
\end{itemize}
\end{corollary}

\begin{proof}
Clearly, ${\,g\in \mathcal{L}^{I}\,}$ is uniquely defined by the
given recurrence relations, since ${\,e\,+\,g\,}$ (\,and hence $\,g\,$\,)
at the right hand side of \eqref{2.8} is evaluated only up to time
$\,t\,=\,n\,\tau\,$\, in virtue of the causality of $\,F\,$.
\newline Furthermore, if
$\,N\leqslant n\,$ and $a\,\leqslant\,\|\,g\,(\,n\tau\,)\,\|\,$
with any $a\,$ satisfying \eqref{2.4}, then
\newline with pertinent $b,c\,$ that satisfy \eqref{2.4}
\begin{equation}\centering\notag 
\begin{aligned}
0\;\leqslant\;(\,\frac{1}{b}\,
&(\,\|\,g\,(\,n \tau\,)\,\|\,-\,a\,)\,)^{\,1/c}\;
\leqslant\;\alpha\,(\,g\,(\,n \tau\,)\,) \\
&=\;\frac{1}{\tau}\;(\;\int_{s < (n+1)\,\tau}\,\alpha\,(\,g\,)\;d\mu(s)
\; - \;\int_{s < n \tau}\,\alpha\,(\,g\,)\;d\mu(s)\;) \\
&=\;\frac{1}{\tau}\;(\;\int_{s < n \tau}\,\alpha\,
(\,F\,[\,e\,+\,g\,]\,)\;d\mu(s)\; - \; 
\int_{s < n \tau}\,\alpha\,(\,g\,)\;d\mu(s)\;) \\
&\;\;\rotatebox{90}{$\rightsquigarrow$}
\text{\small{\;recursion formula \eqref{2.8}}} \\
&\leqslant\;\frac{1}{\tau}\;(\;\int_{s<n\tau}\,\alpha\,(\,e\,+\,g\,)\;d\mu(s)
\;-\;\int_{s<n\tau}\,\alpha\,(\,g\,)\;d\mu(s)\;) \\
&\;\; \rotatebox{90}{$\rightsquigarrow$}
\text{\small{\;since F is passive}} \\
&=\;\frac{1}{\tau}\;\int_{[\,0,\,N\tau\,)}\,\alpha\,(\,e\,+\,g\,)
\,-\,\alpha\,(\,g\,)\;\;\;d\mu\; \\
&\;\; \rotatebox{90}{$\rightsquigarrow$}
\text{\small{\;since $\,e(\,t\,)\equiv 0\,$\; if
$\,t \notin [\,0,\,N\tau)\,$}}\quad.
\end{aligned}
\end{equation}
Thus, estimates \eqref{2.9} applies in the case
$a\,\leqslant\,\|\,g\,(\,n\tau\,)\,\|\,$ and trivially otherwise.
It follows that ${\,\|\,g\,\|\,}$ is uniformly bounded on 
${\,I \smallsetminus [\,0,\,N\tau)\,}$, hence also on $\,I\,$, 
which is to say that $\,g\,$ is stable.
\end{proof}
\vspace{-0.2cm}
\section{Stable DSC processes}\label{S:sec3}

In the following, DSC schemes are represented as paired $\,\alpha$-passive
processes such as dealt with in Theorem ~\ref{T.2.1}, and which hence are
stable.

\begin{lproposition}\label{P.3.1}
For every $\,r,s\in J\,$ and $\,T_{q}:f(t)\,\mapsto\,f(t+q)\,$ denoting the
\emph{shift operator} on $\,\mathcal{L}^{I\;\cup\,J}\,$,
\vspace{-8pt}
\begin{equation}\centering\notag 
\text{if}\quad r+s\;\leqslant\;0\quad\text{and}\quad
\begin{cases}
F_{_{R}}\\[-2pt]
F_{_{C}}
\end{cases}
\!\!\!\!\text{is $\,\alpha$-passive on}\quad
\begin{cases}
\mathcal{L}^{I}\\[-2pt]
\mathcal{L}^{J}
\end{cases}
\!,
\end{equation}
\vspace{-5pt}
\begin{equation}\centering\notag
\text{then}\quad
\begin{cases}
\,T_{r}\circ F_{_{R}}\circ T_{s}\\[-2pt]
\,T_{r}\circ F_{_{C}}\circ T_{s}
\end{cases}
\!\!\!\!\text{is $\,\alpha$-passive on}\quad
\begin{cases}
\mathcal{L}^{J}\\[-2pt]
\mathcal{L}^{I}
\end{cases}
\!.
\end{equation}
\end{lproposition}

\hspace{-.55cm}
The statement is a direct consequence of Definition ~\ref{D.2.3}
\hfill $\Box$

\begin{ldefinition}\label{D.3.1}
\begin{equation}\centering\notag
\text{The}\quad
\begin{cases}
\text{reflection map $\mathcal{R}$}\\[-2pt]
\text{connection map $\mathcal{C}$}
\end{cases}
\text{is called $\,\alpha$-passive,\;iff}\quad
\begin{cases}
F_{_{R}}\\[-2pt]
F_{_{C}}
\end{cases}
\text{is $\,\alpha$-passive}
\end{equation}
in the sense of Definition ~\ref{D.2.3}.
\end{ldefinition}

\begin{ltheorem}\label{T.3.1}
With every time limited excitation, the DSC process generated by
$\,\alpha$-passive reflection and connection maps is uniformly bounded,
hence stable.
\end{ltheorem}

\begin{proof}
It is sufficient to show that every finitely excited DSC process
which is generated by $\,\alpha$-passive $\mathcal{R}\,$ and $\mathcal{C}\,$
can be written as a pair of processes, either of which satisfies
Theorem ~\ref{T.2.1}.

With $\,H\,:\,=\,I\cup J\,=\,{\{k \tau/2\,|\,k\in\mathbb{Z}\}}\,$
and the measure on $\,H\,$ inherited from $\,I\,$ and $\,J\,$ jointly,
the space of all DSC processes is
\begin{equation}\centering\notag 
\begin{aligned}
\mathcal{E}\;&:\,=\;\{\;z=(\,z_{in},\,z_{out}\,)
\in (\mathcal{L}^{2})^{\,H}\quad |
\quad z_{in}((2k+1)\,\tau/2\,)\,=\,z_{in}(\,k\tau\,)\\
\quad &\text{and}\quad
z_{out}(\,k\tau\,)\,=\,z_{out}((2k-1)\,\tau/2\,)\;,\quad
\text{for every $\,k\in\mathbb{Z}$}\;\}\;,
\end{aligned}
\end{equation}
i.e. $\mathcal{E}$ consists of all pairs of functions
${\,z\,=\,(\,z_{in}\,,\,z_{out}\,)\,
\in(\mathcal{L}^2)^{H}\,}$,
such that $\,z_{in}:H\to\mathcal{L}\,$ and
$\,z_{out}:H\to\mathcal{L}\,$, respectively, 'switch' at even and odd
integer multiples of half the time step.
So, there is a natural bijection
\vspace{-.05cm}
\begin{equation}\centering\notag 
\begin{aligned}
\mathcal{E}\;
&\to\;\mathcal{L}^{I}\times\mathcal{L}^{J}\\[-3pt]
z\,=\,(\,z_{in}\,,\,z_{out}\,)\;
&\mapsto\;(\,z_{in}\downharpoonright I\,,\,z_{out}\downharpoonright J\,)\quad,
\end{aligned}
\end{equation}
in virtue of which the first and second components in
$\,\mathcal{E}\,$ can be naturally identified with $\,\mathcal{L}^{I}\,$
and $\,\mathcal{L}^{J}\,$, respectively.

From section ~\ref{S:sec1} (Definition\,\ref{D.1.1}) we recall that for
any incident function
$\,e\,\in\,\mathcal{L}^{I}\times\{\,0\,\}\,\subset\mathcal{E}\,$
supported on a finite interval $\,[\,0\,,\,N\tau\,)\subset I\,$,
the DSC process excited with $\,e\,$ and generated by
$\,\mathcal{R}\,$ and $\,\mathcal{C}\,$ is the well-defined function
$\,z\in\,(\,\mathcal{L}^{2}\,)^{H}\,$ 
such that $\,z\,(\,t\,)\,=\,0\,$ for $\,t\,\leqslant\,0\,$ and recursively
for $\,0\,<\,t\,\in\,H\,$
\vspace{-.2cm}
\begin{equation}\centering\label{3.11}
\begin{aligned}
z\,(\,t\,+\,\frac{\tau}{2}\,)\,=\,
\begin{cases}
(\,z_{in}\,(\,t\,)\,,\;T_{_{-\frac{\tau}{2}}}\,F_{_{R}}\,T_{_{-\frac{\tau}{2}}}
\,[\,e\,+\,z_{in}\,]\,(\,t\,)\,)
&\,\text{if $\,t\,\in\,I\,$}\\
(\,T_{_{-\frac{\tau}{2}}}\,F_{_{C}}\,T_{_{-\frac{\tau}{2}}}\,
[\,z_{out}\,]\,(\,t\,)\,,\;z_{out}(\,t\,)\,)
&\,\text{if $\,t\,\in\,J\,$}\,.
\end{cases}
\end{aligned}
\end{equation}
In fact, $\,z\,$ is uniquely defined by these relations which provide
separate recurrence relations for the two processes
$\,z_{in}\,,\,z_{out}\,$;
for instance for $\,z_{in}\,$ and $\,t\,\in I\,$
\vspace{-4pt}
\begin{equation}\centering\label{3.12}
\begin{aligned}
z_{in}\,(\,t\,+\,\tau\,)\;
&=\; z_{in}\,(\,t^{\sptilde}\!+\,\frac{\tau}{2}\,)
\quad\text{with}\quad t^{\sptilde}:\,=\,t\,+\,\frac{\tau}{2}\,\in\,J &&\,\\
&=\; T_{-\frac{\tau}{2}}F_{_{C}}\,T_{\frac{-\tau}{2}} \hspace{-4.9cm}
\underbrace{[\,z_{out}\,]_{_{\leq t^{\sptilde}}}}_{\hspace{5.8cm}
=\,T_{\frac{\tau}{2}}\;T_{-\frac{\tau}{2}}F_{_{R}}\,T_{-\frac{\tau}{2}}
[e + z_{in}]_{_{\leq t^{\sptilde}\!\!-\frac{\tau}{2}\,=\,t}}\;
\text{by \eqref{3.11}}}
\hspace{-4.9cm}
(\;t^{\sptilde}\,) &&\text{by \eqref{3.11}}\\
&=\; F_{_{C}}\,T_{-\frac{t\tau}{2}}\,F_{_{R}}\,T_{-\frac{\tau}{2}}
[\,e\,+\,z_{in}\,]\,(\,t\,)\quad .\,&&
\end{aligned}
\end{equation}
By causality $\,z_{in}\,$ enters the last expression only up to argument
$\,t\,$, hence \eqref{3.12} is a well defined recursion formula for
$\,z_{in}(\,t\,)\,$, $\,t\in I$.
Moreover, since products of $\,\alpha$-passive operators are
$\alpha$-passive,
${F_{_{C}}\,T_{-\tau/2}\,F_{_{R}}\,T_{-\tau/2}}$ is $\,\alpha$-passive
in virtue of Proposition ~\ref{P.3.1}. Hence Theorem~\ref{T.2.1}
applies to $\,z_{in}\,$ which thus is stable, just as then also is
$\,z_{out}\,=\,F_{_{R}}\,T_{-\tau/2}\,[\,z_{in}\,+\,e\,]\,$.
\end{proof}
\vspace{-0.2cm}
\section{Model equations}\label{S:sec4}

The typical structure of the DSC algorithm has been outlined in the last
section.
Here we deal shortly with the dynamical \emph{model equations}
- ~that the algorithm has ultimately to solve, and from which the connection
and reflection maps $\,\mathcal{C}$, $\mathcal{R}\,$ are
derived.
Most fundamentally, only \emph{total fields\/} $z^p\/$, $z^n$ - ~\emph{not\/},
however, their \emph{incident\/} and \emph{outgoing components separately\/},
enter the DSC model equations.
Accordingly we consider only equations between \emph{total} fields, viz.
sums of incident and reflected fields which are of the following so called
\emph{port} and \emph{node} types
\begin{ldefinition}\label{D.4.1}
For $t\in H\,=\,I\cup J$ and
$z(t)=(\,z_{in}(t)\,,\,z_{out}(t)\,)\in\,
\mathcal{E}\subset(\mathcal{L}^{2}\,)^{H}$,
\vspace{-12pt}
\begin{equation}\centering\label{4.1}
 z_{in}\,(\,t\,)\,+\;z_{out}\,(\,t\,)\quad =\,:\quad
\begin{cases}
z^p\,(\,t\,\,)
\quad\text{if $\,t\,\in\,I\,$}\\[-2pt]
z^n\,(\,t\,\,) 
\quad\text{if $\,t\,\in\,J\,$}\quad.
\end{cases}
\end{equation}
\end{ldefinition}

In general, viz. prescinding from certain restrictions that can be inferred
from a principle of \emph{near field interaction} ~\cite{He3}, the DSC model
equations are of the two types
\begin{equation}\label{4.2}\centering
\mathcal{F}^{n}[z_{+}^{n}][z^{p}]\;\;\;\equiv\,\,0\quad\text{and}\qquad
\mathcal{F}^{p}[z_{+}^{p}][z^{n}]\;\;\;\equiv\,\,0 \quad ,
\end{equation}
with causal functions $\mathcal{F}^{n}$, $\mathcal{F}^{p}$ that respectively
switch on $\,I\,$ and $\,J\,$. The brackets are here short hand for 'back
in time running' sequences of states (along with their \emph{history}\,)
\vspace{-0.18cm}
\begin{equation}\centering\notag
[\,z_{_{\,(\pm)}}\,]\;:\,=\;[\,z\,]_{_{\leq\,t\;(\pm\,\frac{\tau}{2})}}\;=\;
(\,z\,(\,t\;(\pm\,\frac{\tau}{2})\,-\,n\tau)\,)_{_{n\,=\,0,1,2,...}}\;,\quad
t\,\in\,H\;.
\end{equation}
The ${\tau/2}$ time shifts synchronize node and cell boundary switching
in \eqref{4.2} such that the equations hold on all time intervals
$\,{[\,t,\,t+\tau\,)}\,$, ${\,t\in I\,}$ or ${\,t\in\,J\,}$,
respectively, in harmony with the switching conventions ~\eqref{4.1}
for port and node quantities.
Of course, time shifts by ${-\tau/2}$ would also bring about synchronization.
However, the resulting equations would not yield explicit schemes, in general,
or conflict with the causality requirements for $\mathcal{R}$ and 
$\mathcal{C}$.

Inspection of equations \eqref{4.2} shows that $\mathcal{F}^{n}$ affects
only the reflection cycle, while $\mathcal{F}^{p}$ applies only to the
connection cycle. 

The model equations are treated in more detail in ~\cite{He3}. 
In the next section the derivation of DSC model equations is exemplified 
with a heat propagation scheme.
\vspace{-0.3cm}
\section{A heat propagation scheme in non-orthogonal mesh}\label{S:sec5}

The physical interpretation underlying the following application relates
a smoothly varying (viz. in time and space continuously differentiable,
$C^1$-) temperature field $T$,
evaluated as $T^{p}$ at the face centre points and as $T^{n}$
in the nodes of a of non-orthogonal hexahedral mesh, to total states
${ z_{\mu}^{p,n}}$ of a DSC model. 
In fact, we derive the model equations for the connection and reflection
cycles of a DSC heat propagation (diffusion) scheme. 
Since the equations are linear and of dynamical order 0 and 1,
respectively - \nolinebreak as will be seen \nolinebreak - they can be
processed, following the guidelines of the last section.
In the end, we display some computational results of a dispersion test
carried out with this model.

In order to simplify the notation we follow Einstein's convention
to sum up over identical right-hand \nolinebreak (!) sub and
superscripts within all terms where such are present (summation is
not carried out over any index that also appears somewhere at the
left-hand side of a pertinent symbol - \nolinebreak thus, in
${(-1)}^{\kappa} \, a_{\kappa}^{\lambda} \, b_{\lambda} \, \, _{\kappa} c \,$ \,
the sum is made over $\lambda$ but not over $\kappa\,$).

\begin{figure}[!h]\centering
\setlength{\unitlength}{1.cm}
\begin{pspicture}(-1.4,-.5)(20,3.5)\centering
\psset{xunit=.5cm,yunit=.5cm}
\psline[linewidth=0.1mm]{->}(2.5,1.5)(0.0,0.0)
\psline[linewidth=0.1mm]{->}(2.5,1.5)(6.0,2.0)
\psline[linewidth=0.1mm]{->}(2.5,1.5)(3.0,5.0)
\psline[linewidth=0.4mm]{->}(0.0,0.0)(0.0,4.0)
\psline[linewidth=0.4mm]{->}(0.0,0.0)(5.0,0.0)
\psline[linewidth=0.4mm]{->}(6.0,2.0)(6.0,5.0)
\psline[linewidth=0.4mm]{->}(6.0,2.0)(5.0,0.0)
\psline[linewidth=0.4mm]{->}(5.0,0.0)(3.0,3.0)
\psline[linewidth=0.4mm]{->}(3.0,5.0)(0.0,4.0)
\psline[linewidth=0.4mm]{->}(3.0,5.0)(6.0,5.0)
\psline[linewidth=0.4mm]{->}(0.0,4.0)(3.0,3.0)
\psline[linewidth=0.4mm]{->}(6.0,5.0)(3.0,3.0)
\rput(1.8,0.8){$_{_{0}} e$}
\rput(6.0,1.0){$_{_{1}} e$}
\rput(4.0,4.2){$_{_{2}} e$}
\rput(1.5,4.9){$_{_{3}} e$}
\rput(5.0,1.5){$_{_{4}} e$}
\rput(4.4,5.5){$_{_{5}} e$}
\rput(1.5,2.9){$_{_{6}} e$}
\rput(2.5,-.5){$_{_{7}} e$}
\rput(2.4,4.1){$_{_{8}} e$}
\rput(-.6,2.0){$_{_{9}} e$}
\rput(4.0,2.6){$_{_{10}} e$}
\rput(6.7,3.4){$_{_{11}} e$}
\rput(7.0,-0.5){\small{\textnormal{(a)}}}
\psline[linewidth=0.1mm]{->}(15.5,1.5)(13.0,0.0)
\psline[linewidth=0.1mm]{->}(15.5,1.5)(19.0,2.0)
\psline[linewidth=0.1mm]{->}(15.5,1.5)(16.0,5.0)
\psline[linewidth=0.1mm]{->}(13.0,0.0)(13.0,4.0)
\psline[linewidth=0.1mm]{->}(13.0,0.0)(18.0,0.0)
\psline[linewidth=0.1mm]{->}(19.0,2.0)(19.0,5.0)
\psline[linewidth=0.1mm]{->}(19.0,2.0)(18.0,0.0)
\psline[linewidth=0.1mm]{->}(18.0,0.0)(16.0,3.0)
\psline[linewidth=0.1mm]{->}(16.0,5.0)(13.0,4.0)
\psline[linewidth=0.1mm]{->}(16.0,5.0)(19.0,5.0)
\psline[linewidth=0.1mm]{->}(13.0,4.0)(16.0,3.0)
\psline[linewidth=0.1mm]{->}(19.0,5.0)(16.0,3.0)
\psline[showpoints=true,linewidth=0.4mm]{->}(17.375,3.375)(15.0,1.75)
\psline[showpoints=true,linewidth=0.4mm]{->}(14.375,2.625)(18.0,2.50)
\psline[showpoints=true,linewidth=0.4mm]{->}(16.37,.875)(16.0,4.25)
\psline[showpoints=true,
linewidth=0.6mm]{-}(16.1875,2.5625)(16.1875,2.5625) 
\rput(14.55,1.50){$_{_{0}} b$}
\rput(18.5,2.60){$_{_{1}} b$}
\rput(16.6,4.40){$_{_{2}} b$}
\rput(20,-0.5){\small{\textnormal{(b)}}}
\end{pspicture}
\caption{\textsl{Non-orthogonal hexahedral mesh cell. \newline
\textnormal{(a)} Edge vectors. \qquad \qquad 
\textnormal{(b)} Node vectors.}}\label{F:2}
\end{figure}
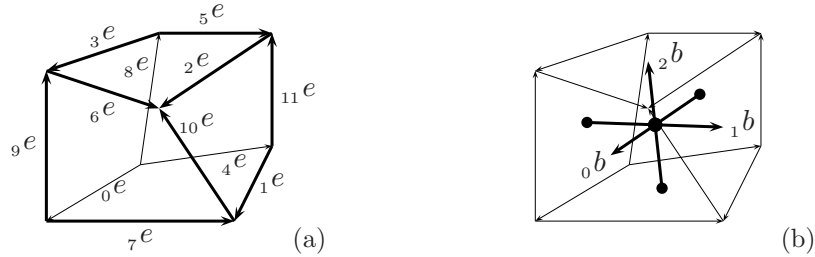

\begin{figure}[!h]\centering
\setlength{\unitlength}{1.cm}
\begin{pspicture}(-2.0,-2)(10,4.0)\centering
\psset{xunit=.5cm,yunit=.5cm}
\psline[linewidth=0.2mm]{->}(0.0,4.0)(0.0,1.0)
\psline[linewidth=0.2mm]{->}(0.0,1.0)(5.0,0.0)
\psline[linewidth=0.2mm]{->}(0.0,4.0)(4.0,5.0)
\psline[linewidth=0.2mm]{->}(4.0,5.0)(5.0,0.0)
\psline[showpoints=true,linewidth=0.6mm]{-}(2.25,2.5)(2.25,2.5) 
\psline[linewidth=0.4mm]{->}(4.5,2.5)(0,2.5)
\psline[linewidth=0.4mm]{->}(2.0,4.5)(2.5,0.5)
\psline[showpoints=true,
linewidth=0.4mm]{->}(4.5,2.5)(8.751,3.400) 
\psline[showpoints=true,linewidth=0.4mm]{->}(0.0,2.5)(-1.5,2.5)
\psline[showpoints=true,linewidth=0.4mm]{->}(2.0,4.5)(1.287,7.352)
\psline[showpoints=true,linewidth=0.4mm]{->}(2.5,0.5)(1.650,-3.751)
\rput(1.0,3.0){$_{_{0}}b$}
\rput(3.0,1.5){$_{_{1}}b$}
\rput(6.6,3.6){$_{_{0}}f$}
\rput(-1.2,3.2){$_{_{1}}f$}
\rput(2.2,6.4){$_{_{2}}f$}
\rput(2.7,-1.4){$_{_{3}}f$}
\rput(5.4,2.0){$_{_{0}}T^{p}$}
\rput(-0.8,1.9){$_{_{1}}T^{p}$}
\rput(1.0,4.7){$_{_{2}}T^{p}$}
\rput(1.6,0.0){$_{_{3}}T^{p}$}
\rput(2.9,3.0){$T^{n}$}
\end{pspicture}
\caption{\textsl{Face vectors and
temperature points (nodal section).}}\label{F:3}
\end{figure}
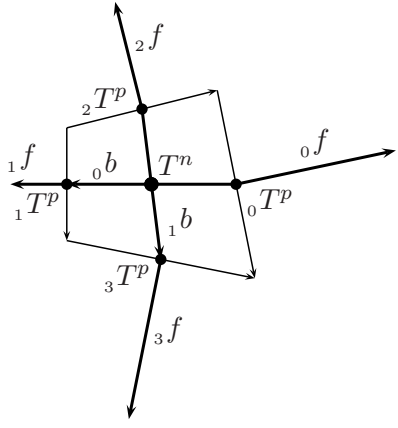

The shape of a hexahedral cell is completely determined by its 12
\emph{edge vectors} ${(_{\nu} e)_{\nu=0,...,11}}\/$.
Also, with the labelling scheme of fig \nolinebreak \ref{F:2}a,
\emph{node
vectors\/} ${(_{\mu} b)_{\mu = 0,1,2}}$ and \emph{face vectors\/}
${(_{\iota} f)_{\iota = 0,...,5}}$ are defined as
\begin{equation}\centering\label{5.1}
\begin{split}
\begin{aligned}
_{\mu} b \quad &:\,= \quad \quad \frac{1}{4}
&& \! \! \sum \nolimits_{\nu = 0}^{3} \, _{_{(4 \mu + \nu)}} e \, \,
&&\mu = 0,1,2 \\
\text{and} \qquad
_{\iota} f \quad &:\,= \quad \, \frac{(-1)^{\iota}}{4}
&& \, ( \, \, _{_{(8 + 2\iota)}} e \,
+ _{_{(9 + 2(\iota + (-1)^{\iota}))}} e \, )
\, \, \land && \\
& &&\; \; \land \, ( \, _{_{(4 + 2\iota)}} e \, 
+ _{_{(5 + 2\iota)}} e \, ) \, \, &&\iota = 0,...,5  \,
\end{aligned}
\end{split}
\end{equation}
(edge vector indices cyclic modulo 12\, and
$\land$ denoting the wedge ('cross') product in $\mathbb{R}^3$).

At every face ${\iota \in \{0,...,5\}}\/$ of a mesh cell and for any 
given $\tau \in \mathbb{R}_{+}\,$ the following time shifted finite
temperature differences in directions ${ _{\mu} b }$ \linebreak 
(\,$\mu = 0,1,2\,$)
form a vector valued function
\begin{equation}\centering\label{5.2}
\begin{split}
_{\iota}\!{\nabla}^{B}T_{\mu}\,(\,t\,,\,\tau\,)\;:\,=\;
\begin{cases}
\,2\,(-1)^{\iota}(\,T^{n}\,_{\mid\,t-\tau/2} 
-\,_{\iota}T^{p}\,_{\mid\,t} \,) \quad
&\text{if $\mu \, = \, [ \iota / 2 ]$} \\[-2pt]
\,(\,\,_{2\mu +1}T^{p}\,-\,_{2\mu}T^{p}\,\,)
\,_{\mid\,t-\tau\,}\quad
&\text{if $\mu\,\neq\,[\iota /2]$}\,
\end{cases}
\end{split}
\end{equation}
($\,[\, x \,]$ denotes the integer part of $x \in \mathbb{R}\,$).
The time increments are chosen to attain technical consistence
with the updating conventions of DSC schemes. They do not destroy
convergence, as easily seen:
In fact, in the centre point of face ${\, \iota \,}$ the vector
${\,_{\iota} \! {\nabla}^{B} T\,}$
approximates in the first order of the time increment ~$\tau$,
and of the linear cell extension, the scalar products of the node
vectors with the temperature gradient ${\, \nabla T \,}$.
Let, precisely, for a fixed centre point on face ${\, \iota \,}$
and $\,\epsilon \in \mathbb{R}_{+}\,$ the \emph{$\epsilon$-scaled cell}
have edge vectors
$\, _{\iota} e\sptilde \, : \, = \, \epsilon \, \, _{\iota} e \,$. 
Let also $\, _{\iota} {\nabla}^{B\sptilde} T_{\mu} \,$ denote
function \eqref{5.2} for the $\epsilon$-scaled cell (with node vectors 
$\, _{\mu} b\sptilde \, = \, \epsilon \, _{\mu} b \;$).
Then at the fixed point holds
\begin{equation}\centering\label{5.3}
\begin{split}
< \, _{\mu} b \, , \, \text{grad($T$)} \, > \, \, \,
= \, \, _{\mu} b \cdot \nabla T \, \,
= \, \, \lim_{\epsilon \to 0} \, \, \lim_{\tau \to 0} \, \,
\frac{1}{\epsilon} \, _{\iota} \! {\nabla}^{B^{\sptilde}} T_{\mu} \, ,
\end{split}
\end{equation}
as immediately follows from the required $C^1$-smoothness of the
temperature field $T$.

To recover, in the same sense and order of approximation, the gradient
${{\nabla} T \,}$ from \eqref{5.2}, observe that for every orthonormal
basis $( _{\nu} u )_{\nu = 0,...,m-1}\,$ of
$V = \mathbb{R}^{m} \, \text{or} \, \, \mathbb{C}^{m}\,$, and for an
arbitrary basis $( _{\mu} b )_{\mu = 0,...,m-1}$ with coordinate
matrix ${\beta_{\nu}^{\mu}} \, = \, {< \, _{\nu} u \, , \, _{\mu} b \, >}$,
the scalar products of any vector $a \in V$ with ${_{\mu} b}$ are
\vspace{-.18cm}
\begin{equation}\centering\label{5.4}
\underbrace{<\,_{\mu} b\,,\,a\,>}_{\qquad=\,: 
\,\,{\alpha}_{\mu}^{B}}\,\,=\,\sum \nolimits_{\nu = 0}^{m-1} \,
\underbrace{<\,_{\mu}b\,,\,_{\nu}u\,>}_{\,\,\,
({\bar{\beta}}_{\mu}^{\nu})\,=\,({\beta}_{\nu}^{\mu})^{^{*}}} \,
\underbrace{<\,_{\nu}u\,,\,a\,>}_{\qquad=\,:\,\,{\alpha}_{\nu}}
\,\,=\,\bar{\beta}_{\mu}^{\nu}\,{\alpha}_{\nu}\;.
\end{equation}
(At the right-hand side, and henceforth,
we follow Einstein's convention to sum up over identical sub and superscripts
within terms where such are present). It follows that
\begin{equation}\centering\label{5.5}
{\alpha}_{\nu}\,=\,
{\gamma}_{\nu}^{\mu} \alpha_{\mu}^{B}\,,
\qquad\text{with}\qquad ({\gamma}_{\nu}^{\mu}) 
\,:\,=\,{({(\beta_{\nu}^{\mu})}^{*})}^{-1}\quad .
\end{equation}

In other words, the scalar products of any vector with the basis vectors
${_{\mu} b\,}$ transform into the coordinates of that vector with
respect to an orthonormal basis ${_{\nu} u \,}$ via multiplication
by matrix ${\gamma = (\beta^*)^{-1}\,}$, where
${\beta_{\nu}^{\mu}}\,\,=\,{<\, _{\nu} u\,,\, _{\mu}b\,>}\,$,
i.e. $\beta$ is the matrix of the coordinate (column) vectors
${_{\mu} b \,}$
with respect to the given ON-basis ${_{\nu} u \,}$, and $\gamma$ is the
adjoint inverse of $\beta$.

This applied to the node vector basis ${ _{\mu} b \,}$ and \eqref{5.3}
yields the approximate temperature gradient at face $\iota$ as
\begin{equation}\centering\label{5.6}
_{\iota} \! \nabla T_{\nu} \quad
= \quad {\gamma}_{\nu}^{\mu}\,\,\, _{\iota}\!{\nabla}^{B} T_{\mu}.
\end{equation}

The heat current \emph{into} the cell through face $\iota$ with
face vector components
${ _{\iota} f^{\nu}}\,
=\,\,{<\,_{\iota} f\,,\,_{\nu} u\,>}\,$,\,
$\nu \in \{ 0,1,2 \}\,$, is then
\begin{equation}\centering\label{5.7}
\begin{aligned}
\qquad _{\iota}J\,\, 
&=\,\,{\lambda}_{H}\,_{\iota}f\,\cdot\,_{\iota}\!\nabla T\,\,
=\,\,\underbrace{{\lambda}_{H}\,_{\iota}f^{\nu}\,\,
{\gamma}_{\nu}^{\mu}}_{\qquad=\,:\,\,_{\iota}s^{\mu}}\,\,
_{\iota}\!{\nabla}^{B}T_{\mu}\,\,
=\,\,_{\iota} s^{\mu}\,\,_{\iota}\!{\nabla}^{B} T_{\mu}\quad,
\end{aligned}
\end{equation}
${{\lambda}_{H}}\,$ denoting the heat conductivity in the cell.

The heat current through every interface is conserved, i.e. between any
two adjacent cells $\zeta$, $\chi$ with the common face labelled $\iota$
in cell $\zeta$ and $\kappa$ in $\chi$ applies
\vspace{-.18cm}
\begin{equation}\centering\label{5.8}
_{\iota}^{^{\zeta}} \! J \quad = \quad - \,\, _{\kappa}^{^{\chi}} \! J \, .
\end{equation}
Also, the nodal temperature change in cell $\zeta$ is
\begin{equation}\centering\label{5.9}
\frac{ d \, {T}^{n}}{dt} \quad 
= \quad \frac{1}{c _{v} \, V } \,
( \, \, S \, + \, \sum \nolimits_{\iota = 0}^{5} \,
_{\iota}^{^{\zeta}} \! J \, \, ) \, ,
\end{equation}
where ${c_{v}}$ denotes the heat capacity (per volume), $V$ the cell
volume, and $S$ any heat source(s) in the cell. \\
We finally introduce quantities
${_{\iota} z _{\mu}^{p,n}}$ \, ($\iota \, = \, 0,...,5;$ \, $\mu \, = \,
0,1,2 \,$), which still smoothly vary in time with the
temperature $T$ (and that are hence not yet DSC states, but will be later 
updated as such)
\begin{equation}\centering\label{5.10}
\begin{split}
_{\iota} z_{\mu}^{n}\,(\,t\,)\quad :\,=\quad
\begin{cases}
\,\,2\,(-1)^{\iota}\,\,T^{n}\,_{\mid\,t}\qquad
&\text{if $\mu\,=\,[\iota /2]$}\,\,\\
\,\,(\,_{2\mu +1} T^{p} 
-\,_{2\mu} T^{p}\,)_{\mid\,t-\tau/2}\qquad
&\text{else} \,
\end{cases}
\end{split}
\end{equation}
and
\vspace{-0.25cm}
\begin{equation}\centering\label{5.11}
\begin{split}
_{\iota}z_{\mu}^{p}\,(\,t\,)\quad :\,=\quad
\begin{cases}
\,\,2\,(-1)^{\iota}\,\,_{\iota}T^{p}\,_{\mid\,t}\qquad
\qquad\;&\qquad\!\!\text{if $\mu\,=\,[\iota /2]$}\;\\
\,\,\,_{\iota} z_{\mu}^{n}\,(\,t\,-\,\tau/2\,)\qquad
\qquad\;&\qquad\!\!\text{else}\;.
\end{cases}
\end{split}
\end{equation}
From (\ref{5.2}, \ref{5.7}) follows
\vspace{-.18cm}
\begin{equation}\centering\label{5.12}
\begin{split}
_{\iota} J\, _{\mid \, t + \tau / 2} \quad
&= \quad \, _{\iota} s^{\mu}\, ( \, _{\iota} z_{\mu}^{n}\, _{\mid \, t} \,
- \, 2 \, {(-1)}^{\iota} {\delta}_{\mu}^{[\iota/2]} \, \,
_{\iota} T^{p} \,_{\mid \, t + \tau/2}\, ) \\
&= \quad _{\iota} s^{\mu}\, ( \, _{\iota} z_{\mu}^{n}\, _{\mid \, t} \,
- \, {\delta}_{\mu}^{[\iota/2]} \, \, 
_{\iota} z_{\mu}^{p}\, _{\mid \, t + \tau/2} \, )
\quad .
\end{split}
\end{equation}
This, together with \eqref{5.8}, \eqref{5.11}, and continuity of the
temperature at the interface,
$\,_{\iota}^{^{\zeta}} T\,^{p}\,=\,_{\kappa}^{^{\chi}} T\,^{p}\,$,
implies
\begin{equation}\centering\label{5.13}
\begin{split}
_{\iota}^{^{\zeta}} z\,_{\mu}^{p}\,_{\mid\,t+\tau/2}\,\,=\,
\begin{cases}
\,\frac{\,_{\iota}^{^{\zeta}} s\,^{\mu}\,\, 
_{\iota}^{^{\zeta}} z\,_{\mu}^{n}\,_{\mid\,t} \,
+\,_{\kappa}^{^{\chi}} s\,^{\nu}\,\,\,
_{\kappa}^{^{\chi}} z \,_{\nu}^{n}\,_{\mid\,t}} 
{_{\iota}^{^{\zeta}} s\,^{[\iota/2]}
+ \, ( -1 )^{\iota +\kappa}\,_{\kappa}^{^{\chi}} s\,^{[\kappa /2]}}
\qquad\,\,
&\text{if $\mu = [\iota/2]$}\vspace{0.2cm}\\
\,\,_{\iota}^{^{\zeta}} z\,_{\mu}^{n}\,_{\mid\,t}
\qquad\,\,
&\text{else}\,,
\end{cases}
\end{split}
\end{equation}
which form a complete set of recurrence relations for ${z^{p}}$
(\,given ${z^{n}}\/$). So they can be taken as model equations for
the connection cycle of a DSC algorithm.

Equations \eqref{5.9} discretely integrated in a time balanced form with
increment ~$\tau$ yield
\vspace{-.18cm}
\begin{equation}\centering\label{5.14}
\begin{split}
T^{n} \, _{\mid \, t \, + \tau/2 } \quad
&= \quad
T^{n} \, _{\mid \, t \, - \tau/2 } \,
+ \frac{\tau}{c_{v} \, V} \, ( \, \, S \, + \, \, \sum \nolimits_{\iota = 0}^{5}
\, _{\iota} J \, _{\mid t } \, \, ),
\end{split}
\end{equation}
i.e., with (\ref{5.10}, \ref{5.11}, \ref{5.12}), the recurrence relations
\vspace{-.10cm}
\begin{equation}
\begin{split}\centering\label{5.15}
_{\iota}z\,_{\mu}^{n}\,_{\mid\,t+\tau/2}\,=\,
\begin{cases}
\,_{\iota} z\,_{[\iota /2]}^{n}\,_{\mid\,t-\tau /2}\,
+\frac{{(-1)}^{\iota}\,\tau}{2\,c_{v}\,V}\,\,\{\,\,S\,\,\,\,+\\
\quad\quad +\,\,\sum\nolimits_{\iota = 0}^{5}\,
_{\iota} s^{\nu}\,(\,_{\iota} z\,_{\nu}^{n}\,_{\mid\,t-\tau/2}\,
-\,{\delta}_{\nu}^{[\iota/2]}\, 
_{\iota} z\,_{\nu}^{p}\,_{\mid\,t}\,)\,\} 
&\text{if $\mu = [ \iota / 2 ]$} \vspace{.2cm}\\
\,-\,\frac{1}{2}\,(\,_{2\mu +1} z\,_{\mu}^{p}\,+
_{2\mu} z\,_{\mu}^{p}\,)\,_{\mid\,t}&\text{else}\,,
\end{cases}
\end{split}
\end{equation}
which provide a complete set of model equations for the reflection cycle
of a DSC algorithm. Note that the first line, modulo the factor
$2\,(-1)^{\iota}\,$, always updates the nodal temperature. Of course,
this has to be carried out only once per cell and iteration cycle.
In this - ~typical ~- example the dual state space concept of DSC
(\,needed by Johns' cycle\,) creates a redundancy, which can yet
be exploited for process parallelization within either parts of the
connection-reflection cycle.

Equations (\ref{5.13}, \ref{5.15}) can be directly taken as updating
relations for total quantities of a DSC scheme. Alternatively, they 
may be further processed in deriving reflection and connection maps,
along with stability bounds for the time step.
The proceeding is canonical and amounts in essence to a straightforward
transcription of the model equations along the lines of ~\cite{He3},
section 4 (Theorem 1 and corollaries).

It is quite easy to couple this heat conduction model - ~within one and
the same mesh ~- to a Maxwell field TLM model in the non-orthogonal
setting \cite{He1}. In fact, with the node vector definition in \cite{He1},
the total node voltages are just the scalar products of $_{\mu} b$ with the
electric field, hence the dielectric losses and heat sources per cell are
\begin{equation}\centering\label{5.16}
S \quad = \quad \frac{1}{2} \,
\sigma \, V \, E^{\nu} \, \overline{E_{\nu}} \quad
= \quad \frac{1}{2} \, \sigma \, V \, \sum \nolimits_{\nu} \, \left|
{\gamma}_{\nu}^{\mu} \, U_{\mu}^{n} \right|^{2} \; ,
\end{equation}
for a frequency domain (complex) TLM algorithm, cf.\cite{He2};
$\sigma=2\pi f\,\epsilon\,\text{tan}(\delta)$
denotes the effective loss current conductivity at frequency $f$
in a mesh cell of absolute permittivity $\epsilon$ and dielectric loss
factor $\text{tan}( \delta )\,$; $\gamma \, = \, ( \beta^{*} )^{-1}\,$
as in \nolinebreak \eqref{5.5}. In \textsc{Spinner}'s Maxwell field
solver the model couples, in addition, to magnetic and skin effect losses.

\begin{figure}[!h]\centering
\setlength{\unitlength}{1.cm}
\begin{pspicture}(0.0,0.0)(13.0,4.5)\centering
\psset{xunit=1.0cm,yunit=1.0cm}
\rput(1.7,1.67){\includegraphics[scale=0.3150,clip=0]{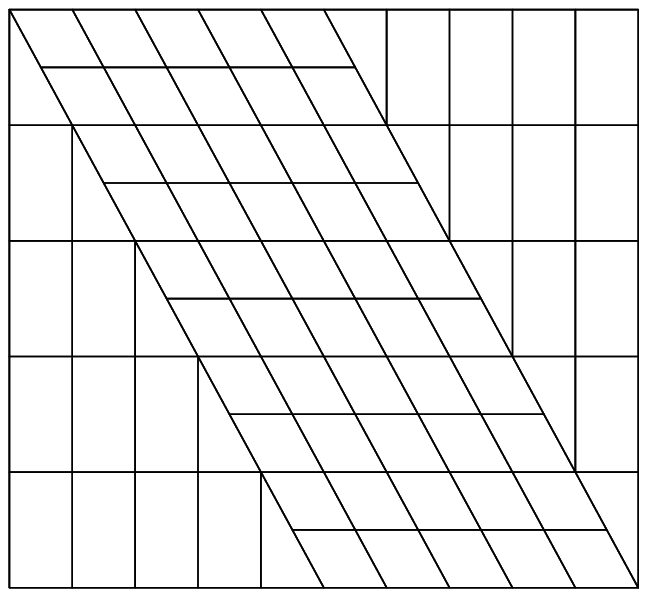}}
\rput(0.25,0.75){\small (a)}
\rput(5.8,1.45){\includegraphics[scale=0.3000,clip=0]{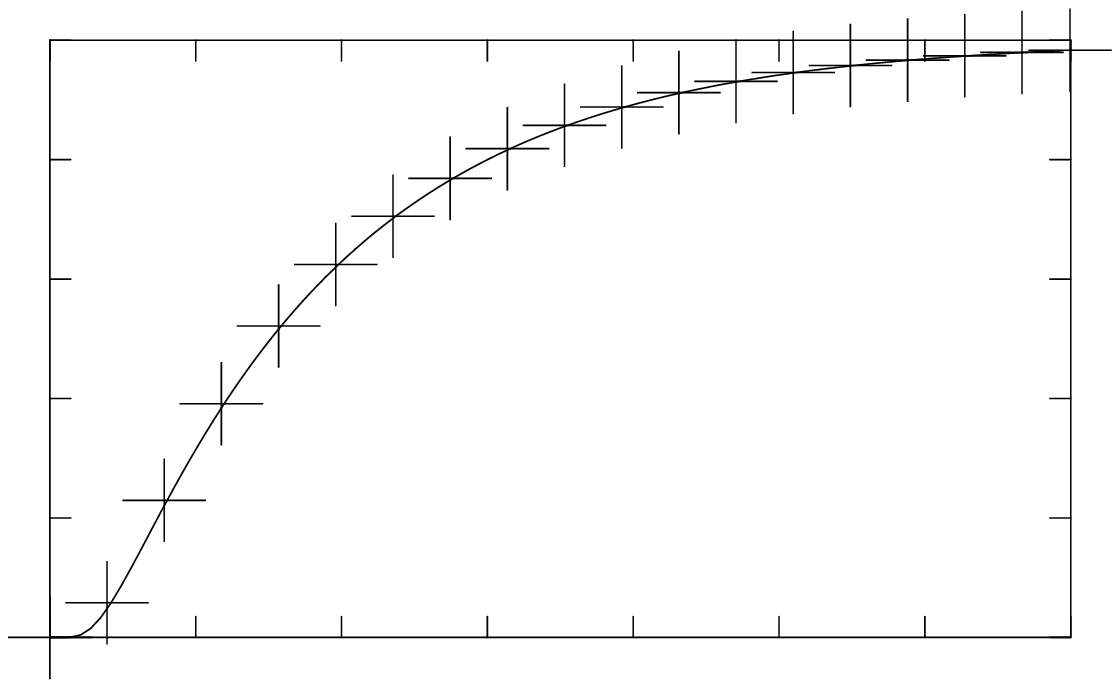}}
\rput(4.0,0.75){\small (b)}
\rput(10.6,1.45){\includegraphics[scale=0.3000,clip=0]{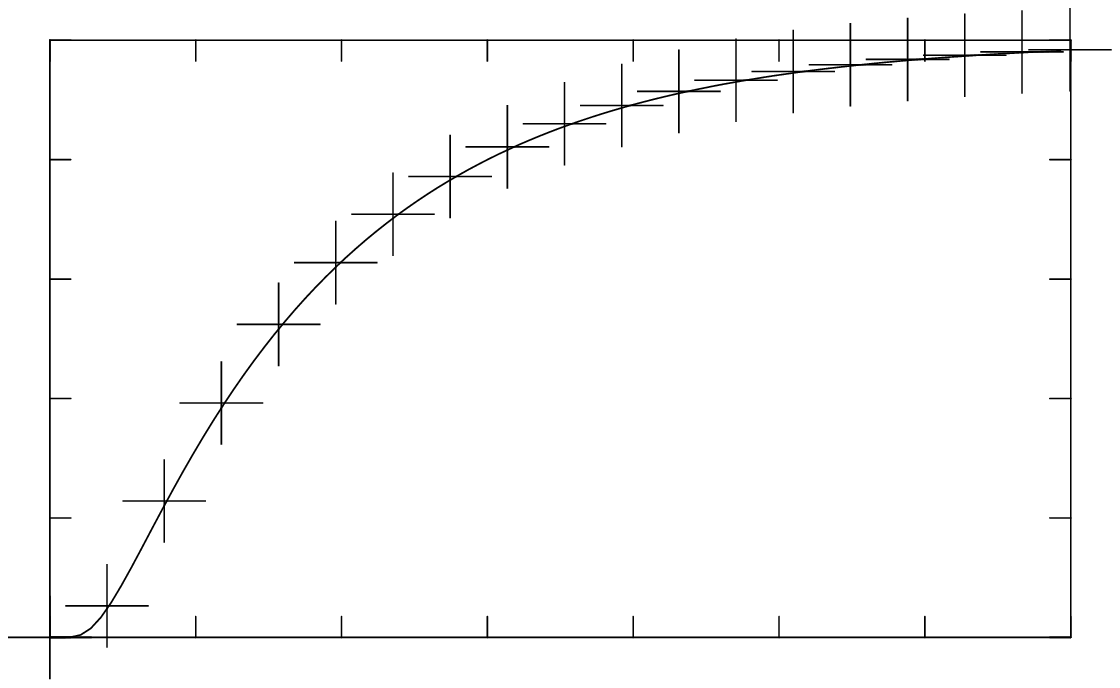}}
\rput(8.8,0.75){\small (c)}
\end{pspicture}
\vspace{-0.8cm}
\caption{\textsl{Transverse heat conduction over a square mesh using
non-orthogonal cells. A Heaviside temperature step is imposed on one
side and the transient temperature computed at the opposite side,
assuming adiabatic boundary conditions on all but the heated sides.
\newline
DSC results \textnormal{(+)} are plotted over analytical solution
(curve).
\newline
\textnormal{(a)} The mesh. \, 
\textnormal{(b)} Horizontal \quad
\textnormal{(c)} vertical propagation.}}\label{F:4}
\end{figure}

Fig \nolinebreak \ref{F:4} displays the result of a dispersion test,
computed in a square mesh using non-orthogonal cells.
It turns out that the heat conduction properties of the mesh are highly
insensitive to cell shape and orientation (as of course should be the case).

\vspace{-0.2cm}
\section{Conclusions}\label{C:sec6}

DSC schemes extend the Transmission Line Matrix numerical method in preserving
its essential assets. So, they remain stable under quite general circumstances
(made tangible with the notion of $\alpha$-passivity) and are amenable to
deflection techniques which generalize the stub-loading feature of TLM.
DSC schemes offer enhanced modeling potentiality, and a wide field of future
research.
\vspace{-0.2cm}

\vspace{3pt} 
\hspace{-16pt}{
\Small{\address{\textsc{Spinner RF Lab, Aiblinger Str.30,
DE-83620 Westerham, Germany}}}

\hspace{-16pt}{
\Small{Email address:\; \email{s.hein@spinner.de}}}
\end{document}